\documentclass[letterpaper,11pt]{article}

\usepackage[latin1]{inputenc}
\usepackage{amsmath, amsfonts, amssymb, amsthm, multirow, numprint, epsfig, bbm}

\newtheorem{thm}{Theorem}%[section]

\newtheorem{lemma}[thm]{Lemma}
\newtheorem{prop}[thm]{Proposition}

\newtheorem{question}[thm]{Question}

\newcommand{\beq}[1]{\begin{equation}\label{#1}}
\newcommand{\enq}[0]{\end{equation}}

\newcommand{\bn}[0]{\bigskip\noindent}
\newcommand{\mn}[0]{\medskip\noindent}
\newcommand{\nin}[0]{\noindent}

\newcommand{\sub}[0]{\subseteq}
\newcommand{\sm}[0]{\setminus}
\renewcommand{\dots}[0]{,\ldots,}

\newcommand{\ov}[0]{\overline}

\newcommand{\A}[0]{{\cal A}}
\newcommand{\B}[0]{{\cal B}}
\newcommand{\cee}[0]{{\cal C}}

\newcommand{\f}[0]{{\cal F}}

\newcommand{\h}[0]{{\cal H}}

\newcommand{\K}[0]{{\cal K}}

\newcommand{\sss}[0]{{\cal S}}
\newcommand{\T}[0]{{\cal T}}

\newcommand{\ra}[0]{\rightarrow}

\newcommand{\mm}[0]{m}

\newcommand{\0}[0]{\emptyset}

\newcommand{\qqed}[0]{\begin{flushright} \rule{2mm}{3mm} \end{flushright}}
\def\qqqed{\null\nobreak\hfill\hbox{\rule{2mm}{3mm} }\par\smallskip}

\newcommand{\C}[2]{{{#1}\choose{{#2}}}}
\newcommand{\Cc}[0]{\tbinom}

\newcommand{\gb}[0]{\beta }
\newcommand{\gc}[0]{\gamma }

\newcommand{\gD}[0]{\Delta }
\newcommand{\gG}[0]{\Gamma }

\newcommand{\gl}[0]{\lambda }
\newcommand{\gL}[0]{\Lambda}

\newcommand{\gO}[0]{\Omega}

\newcommand{\gz}[0]{\zeta}
\newcommand{\eps}[0]{\varepsilon }
\newcommand{\vt}[0]{\vartheta}

\newcommand{\nex}[0]{_{\ov{x}}}

\newcommand{\ekr}{Erd\H{o}s-Ko-Rado}

\newcommand{\comments}[1]{}
\usepackage[usenames]{color}

\begin{document}

\renewcommand{\thefootnote}{\fnsymbol{footnote}}

\title{On ``stability" in the \ekr{} Theorem\footnotemark}
\author{Pat Devlin, Jeff Kahn}
\date{}
\footnotetext{$^*$ Supported by NSF grant DMS1201337.}

\maketitle
\begin{abstract}
Denote by $K_p(n,k)$ the random subgraph of the usual Kneser graph $K(n,k)$
%(whose vertices are the $k$-subsets of $[n]$)
in
which edges appear independently, each with probability $p$.
%In this note,
Answering a question of Bollob\'as, Narayanan, and Raigorodskii,
we show that there is a fixed $p<1$ such that a.s.
(i.e., with probability tending to 1 as $k\ra\infty$)
the maximum independent sets of $K_p(2k+1, k)$ are precisely 
%``stars"
the sets $\{A\in V(K(2k+1,k)): x\in A\}$ ($x\in [2k+1]$).

We also complete the determination of the order of magnitude
of the ``threshold" for the above property for general $k$ and $n\geq 2k+2$.
This is new for $k\sim n/2$, while
for smaller $k $ it is a recent result of Das and Tran.
\end{abstract}

\section{Introduction}\label{Intro}

The broad context of this paper is
an effort, which has been one of the most interesting and successful
combinatorial trends of the last couple decades, to understand
how far some of the subject's classical results remain true
in a random setting.
Since several nice accounts of these developments are available,
%(and since short papers should have short introductions), 
we will not attempt a review
(see, for example, the survey
\cite{Rodl-Schacht}
or \cite{BNR,BBN}
for discussions
closer to present concerns)
and mainly confine ourselves to the problem at hand.

\medskip
Recall
that, for integers $0 < k < n/2$, the {\em Kneser graph}, $K(n,k)$
has vertices the $k$-subsets of $[n]:=\{1,2\dots n\}$,
with
%vertex set $\C{[n]}{k}$
%(where $[n]=\{1,2\dots n\}$ and $\C{X}{k} = \{\mbox{$k$-subsets of $X$}\}$) with
two vertices adjacent if and only if they are disjoint sets.
In what follows we set $\K=\C{[n]}{k}$ (the vertex set of $K(n,k)$).
A {\em star} is one of the sets $\K_x:=\{A:x\in A\}$ ($x\in [n]$).
We also set $M =\C{n-1}{k-1}$
(the size of a star) and write $\cee$ for the collection of $M$-subsets of $\K$
that are {\em not} stars.

In Kneser-graph terms the classical Erd\H{o}s-Ko-Rado Theorem \cite{EKR}
says that for $k<n/2$,
the independence number of $K(n,k)$ is
%equal to
$M$ and,
moreover, the only independent sets of this size are the stars.

Say a spanning subgraph $H$ of $K(n,k)$ has the {\em EKR property} or {\em is EKR}
if each of its largest independent sets is a star.
We are interested in this property for $H = K_p(n,k)$, the random subgraph of $K(n,k)$
in which edges appear independently, each with probability $p$.
%(So we are talking about {\em percolation} on $K(n,k)$, but will not use this language.)
In particular we are interested in a question
suggested and first studied by
Bollob\'{a}s, Narayanan, and Raigorodskii \cite{BNR}, {\em viz.}
\begin{question}\label{BNRQ}
For what $p=p(n,k)$ is $K_p(n,k)$ likely to be EKR?
\end{question}

Formally,
we would like to estimate
the ``threshold," $p_c=p_c(n,k)$, which we define to be the unique $p$ satisfying
\beq{prob}
\Pr(\mbox{$K_p(n,k)$ is EKR}) = 1/2
\enq
(which does turn out to be a threshold in the
original Erd\H{o}s-R\'enyi sense).
Ideally (or nearly so) one hopes to identify some $p_0$, necessarily close to $p_c$,
such that for fixed $\eps>0$,
$K_p(n,k)$ is a.s. EKR if $p> (1+\eps)p_0$ and a.s. {\em not} EKR if
$p< (1-\eps)p_0$.
(As usual a property holds {\em almost surely} (a.s.) if its
probability tends to 1 as the relevant parameter---here $n$---tends to infinity).

Successively stronger results
(some of this ``ideal" type, some less precise) have been achieved
%for  various ranges of $n$ and $k$
by the aforementioned Bollob\'as {\em et al.}
\cite{BNR} and then by
Balogh, Bollob\'{a}s, and Narayanan \cite{BBN} and
Das and Tran \cite{Das}.
Here we briefly discuss only \cite{Das}, which subsumes the others.

A natural guess is that the value of $p_c$ is driven by the need
to avoid independence of any  $\f\in \cee$ that, for some $x$,
satisfies $|\f\sm\K_x|=1$.
This turns out to suggest that $p_c=p_c(n,k)$ should be asymptotic to
\[
p_0= p_0(n,k) :=\left\{\begin{array}{cl}
\Cc{n-k-1}{k-1}^{-1}\log(n\Cc{n-1}{k}) &\mbox{if $n\geq 2k+2$,}\\
3/4& \mbox{if $n=2k+1$}
\end{array}\right.
\]
(where, here and throughout, $\log $ is $\ln$);
namely,
\cite{Das} shows (strictly speaking only for $n\geq 2k+2$)
that for $p<(1-\eps)p_0$ (with $\eps >0$ fixed), $K_p(n,k)$ a.s. does contain
independent $\f$'s as above (implying $p_c>(1-o(1))p_0$), while it is easy
to see that for $p>(1+\eps)p_0$ it a.s. does not.
(Note $n=2k+1$ is not really special here: the form of $p_0$ changes because
we lose the approximation of $1-p$ by $e^{-p}$.)

In fact, Das and Tran show that, for some specified constant $C$,
$p_c$ is indeed asymptotic
to $p_0$ if $k< n/(3C)$, and is less than
$Cn p_0/(n-2k)$
whenever $n\geq 2k+2$, whence
$p_c =O(p_0)$ if $k < (1/2-\gO(1))n$
(where the first implied constant depends on the second).
Of course the estimate becomes less satisfactory as $k/n\ra 1/2$,
and in particular gives nothing for $n \in \{2k+2,2k+1\}$.
On the other hand, both
\cite{BNR} and \cite{BBN} suggest
that $n=2k+1$ is the most interesting case of the problem
and ask whether one can
at least show that $K_p(2k+1,k)$ is a.s. EKR for some $p$ bounded away from 1.
Here we prove such a result and also show that $p_c =O(p_0)$ remains true
for general $n$ and $k$.

\begin{thm}\label{MT}
There is a fixed $p<1$ such that (for every k)
$K_p(2k+1,k)$ is a.s. EKR.
%{\rm (b)}
There is a fixed $C$ such that for every n and k,
$K_p(n,k)$ is a.s. EKR for $p> Cp_0(n,k)$.
\end{thm}

\nin
Again, one expects that $p_c\sim p_0$ in all cases and in
particular that, as suggested in
\cite{BBN}, $p_c(2k+1,k) \ra 3/4$;
but we do not come close to these asymptotics and make no attempt to
squeeze the best possible $\eps$ and/or $C$ from our arguments.

\medskip
It may be worth (briefly) comparing the present question with a similar
one, introduced earlier by
Balogh, Bohman and Mubayi \cite{BBM}, 
in which one considers a random {\em induced} subgraph of $K(n,k)$.
Thus, one specifies only $\h =\K_p$ (the random set in which each $A\in \K$
is present with probability $p$, independent of other choices) and asks when the
subgraph induced by $\h$ has the {\em EKR property}, now meaning
that each largest independent set (that is,
intersecting subfamily of $\h$) is a {\em star} $\h_x=\{A\in \h:x\in A\}$ for some $x$.

For $n=2k+1$ the situation here is similar to the one above:
EKR should hold (a.s.) for any fixed $p>3/4$, but even proving this for 
$p> 1-\eps$ with a fixed $\eps >0$---a problem suggested in \cite{BBM}---does not
seem easy; such a proof was given in 
\cite{EKRII} (using methods unrelated to those employed here).

But the resemblance may be
superficial, and in fact the induced problem seems 
considerably subtler than the one considered here (as should probably be expected, e.g. since
(i) the size of a largest star is itself a moving target and (ii) the most likely violators
of EKR are not always families that are close to stars).
See \cite{EKRI} for a guess as to what ought to be true here
and \cite{BDDLS,EKRI} for what's known at this time.

\medskip
The rest of the paper is devoted to the proof of
Theorem~\ref{MT}.  A single argument will suffice for both assertions, though,
as noted below, not all of what we do is needed for $n=2k+1$.
%somewhat less would be required if we were only dealing with $n=2k+1$.

\section{Proof}

{\em Notation.}
From now on we take $n=2k+c$ and
write
$V$ for $[n]$ (so $\K=\C{V}{k}$).
For $\h\sub \K$,
we let $\h_x=\{A\in \h:x\in A\}$,
$\h\nex = \h\sm\h_x$ ($x\in V$) and $\gD_\h=\max\{|\h_x|:x\in V\}$.
As usual $|\h_x|$ is the {\em degree} of $x$ in $\h$.
We use $M$ and $\cee$ as above and set
$N=\C{n-1}{k}$.
For $\f\in\cee$ we set
%{\bf OR}  Set $\cee = \C{\K}{M}\sm \{\K_x:x\in [n]\}$.
$a_\f=M-\gD_\f$ and $e(\f)=|\{\{A,B\}:A,B\in \f, A\cap B=\0\}|$
(the number of Kneser edges in $\f$).

\bigskip
In view of \cite{Das} we may assume
\beq{k}
k>6c.
\enq
We also assume henceforth that
\beq{p}
p>\left\{\begin{array}{ll}
1-\eps&\mbox{if $c=1$,}\\
Cp_0(n,k)&\mbox{if $c\geq 2$}
\end{array}\right.
\enq
for suitable fixed $C$ and $\eps>0$
(namely, ones that support our arguments)
and want to show that then
\[
\Pr(\mbox{some $\f\in\cee$ is independent in $K_p(n,k)$})=o(1).
\]
Perhaps surprisingly, this is given by a straight union bound; that is,
there are lower bounds on the sizes of the various $e(\f)$'s
that imply
(with $\f$ running over $\cee$)
\beq{toshow}
%\sum\Pr(\mbox{$\f$ independent})
\sum (1-p)^{e(\f)}=o(1).
%~~( = \sum (1-p)^{|E(\f)|})~~=o(1).
\enq
This contrasts with (e.g.) \cite{EKRII},
where a naive union bound gives nothing.

\medskip
The rest of our discussion is devoted to the proof of \eqref{toshow},
and we assume
from now on that $\f\in\cee$.
Notice that we always have
\beq{abd}
a_\f/N\leq k/n
\enq
(since the trivial $\gD_\f \geq k|\f|/n = kM/n$ gives
$a_\f\leq (1-k/n)M = kN/n$).

\medskip
The next assertion is the main point.

\begin{lemma}\label{ML}
There is a fixed $\vt>0$ such that for any $\f\in \cee$,
\beq{e(f)}
e(\f) > \vt k^{-1}\Cc{n-k-1}{k-1} a_\f\log(N/a_\f).
\enq
\end{lemma}

%\medskip
We first observe that this easily gives \eqref{toshow}.
Noting that (for any $a$) the number of $\f$'s with
$a_\f =a$ is at most $n\C{M}{a}\C{N}{a}$
(choose a maximum degree vertex $x$ of $\f$ and then the
$a$-sets $\K_x\sm\f\sub \K_x$ and $\f\nex\sub\K\nex$), we find that,
with $\vt$ as in Lemma~\ref{ML}, the sum in \eqref{toshow} is then
less than
\beq{eFbound}
n\sum\left\{\Cc{M}{a}\Cc{N}{a}\exp[-\xi\vt k^{-1}\Cc{n-k-1}{k-1}
a\log (N/a)]:0<a\leq kN/n\right\},
\enq
where 
$
\xi=\log(1/\eps)$ if $n=2k+1$
and otherwise $\xi=p$.
%\[\xi=\left\{\begin{array}{ll}
%\log(1/\eps)&\mbox{if $n=2k+1$,}\\
%p&\mbox{otherwise.}
%\end{array}\right.\]
We may bound the summand using $\C{M}{a}\C{N}{a}< \exp[2a\log (eN/a)]$ and 
\[
\xi\vt k^{-1}\Cc{n-k-1}{k-1} \geq \left\{\begin{array}{ll}
\vt \log(1/\eps)
&\mbox{if $n=2k+1$,}
\\
C \vt k^{-1}\log(n\Cc{n-1}{k})&\mbox{otherwise,}
\end{array}\right.
\]
and the expression in \eqref{eFbound} is then
easily seen to be small if (say) $\eps < e^{-5/\vt}$ or
$C>4/\vt$ (for $n=2k+1$ and $n\geq 2k+2$ respectively).
\qqed

\medskip
The proof of Lemma~\ref{ML} divides into three regimes,
depending on $a_\f$.
The first of these---$a_\f$ not too small---is handled
as in \cite{Das}, from which we recall only what
we need
(see their Theorem 1.2):

\begin{thm}\label{thmDT}
There is a fixed K such that for any $2\leq k<n/2$:
if $\f\in \cee$ satisfies $a_\f >K\gz \frac{n}{c} M$ with
$\gz \leq \frac{c}{(10K)^2n}$, then
$e(\f) > \gz M \C{n-k-1}{k-1}$.
\end{thm}

\mn
It will be convenient to assume (as we may) that $K\geq 1$.
Theorem~\ref{thmDT} gives \eqref{e(f)}
for any $\f$ satisfying $a_\f > M/(100K)$
(with $\vt $ something like $.01K^{-2}$),
so we assume from now on that this is not the case.

\medskip
For smaller values we need to say a little about graphs belonging
to the ``Johnson scheme" (e.g. \cite{MacW-S}).
For positive integers $k\leq m$
we use $J_i(m,k)$ for the graph on
$V_{m,k}:=\C{[m]}{k}$ with $A,B$ adjacent ($A\sim _iB$) iff
$|A\Delta B|=2i$.
Here we take $m=n-1$
and will be interested in $i\in \{1,c\}$.
Uniform measure on $V_{m,k}$ will be denoted $\mu_k$.

We use $\gb_i(\A)$ for the size of the {\em edge boundary} of $\A\sub V_{\mm,k}$ in $J_i(m,k)$;
that is,
\[
\gb_i(\A) = |\{\{A,A'\}: A\in \A, ~A'\in (V_{\mm,k}\sm \A),~ A\sim_i A'\}|.
\]
The following
lower bounds on $\gb_c$ and $\gb_1$ will suffice for our purposes.

\medskip
For $\gb_c$ we use a standard version of the eigenvalue-expansion connection
due to Alon and Milman \cite{Alon-Milman}
(see e.g. \cite[Theorem 9.2.1]{AS}),
which (here) says that for any $\A\sub V_{m,k}$,
\beq{A-M}
\gb_c(\A) \geq \gl |\A|(1-\mu_k(\A)),
\enq
with $\gl$ the smallest positive eigenvalue of
the Laplacian of $J_c(m,k)$
(the matrix $DI_N-A$, where $D= \C{k}{c}\C{m-k}{c}$ and $A$
are
the degree and adjacency matrix of $J_c(m,k)$).
We assert that (assuming \eqref{k})
\beq{lambda}
\lambda = \frac{m}{k}\C{k}{c} \C{m-k-1}{c-1}.
\enq
{\em Proof.}
(This ought to be known, but we couldn't find a reference.)
The eigenvalues of $A$ are (again, see e.g. \cite{MacW-S})
%\[\mbox{$\gl_j:=\sum_{i=0}^c(-1)^i\Cc{j}{i}\Cc{k-j}{c-i}\Cc{m-k-j}{c-i} ~~~~ j= 0\dots k$}\]
\[
\mbox{$\gl_j:=\sum_{i=0}^c(-1)^i S^j_i,  ~~~~ j= 0\dots k$},
\]
where $S^j_i:=
\Cc{j}{i}\Cc{k-j}{c-i}\Cc{m-k-j}{c-i} $.
In particular,
$
\gl_{0} = \C{k}{c}\C{m-k}{c}$, 
\beq{lambda1}
\gl_{1} = \Cc{k-1}{c}\Cc{m-k-1}{c}-\Cc{k-1}{c-1}\Cc{m-k-1}{c-1} 
=\Cc{k-1}{c}\Cc{m-k-1}{c}\tfrac{km-k^2-cm}{km-k^2-cm+c^2}
\enq
and $\gl_0-\gl_1 =\gl$ (the value in \eqref{lambda}), so we just need to 
show that $\gl_j\leq \gl_1$ for $j\geq 2$.
In fact it is enough to show that 
\beq{Sji}
\mbox{$S^j_i\leq \gl_1~$ whenever $j\geq 2$,} 
\enq
since
log-concavity of 
the sequences $\left(\C{a}{\ell}\right)_\ell$ implies log-concavity
of 
$(S^j_i)_i$ and thus
$ 
\gl_j\leq \max_iS^j_i .
$ 

Routine manipulations (using the expression for $\gl_1$ in \eqref{lambda1}) give
\[   
S^j_i/ \gl_1
=\dfrac{km-k^2}{km-k^2-cm}\dfrac{\C{c}{i}\C{k-c}{j-i}}{\C{k}{j}}\dfrac{\C{m-k-c}{j-i} \C{c}{i}}{\C{m-k}{j}}\dfrac{1}{\C{j}{i}} \\
\leq \dfrac{1}{\C{j}{i}}
%{\C{j}{i}}^{-1} 
\dfrac{km-k^2}{km-k^2-cm}.
\]   
For $0<i<j$, the r.h.s. is less than 1 since $km-k^2 < 2(km-k^2-cm)$, as follows easily
from \eqref{k}.  On the other hand, it is easy to see (using \eqref{k}) that each of 
$S^2_0$ and $S^2_2$ is less than $\gl_1$, which gives \eqref{Sji} 
for $i\in \{0,j\}$, since $S^j_0$ and $S^j_j$ are decreasing in $j$.

\qqed

For $\gb_1$ we
use an instance of a result of Lee and Yau \cite{Lee-Yau} (estimating the log-Sobolev constant
for $J_1(m,k)$):  there is a fixed $\gc>0$ such that, for any
$k$ as in \eqref{k}
and $\A\sub \C{[\mm]}{k} $,
\beq{betalb}
\gb_1(\A) > \gc \mm |\A|
\log ( 1/\mu_k(\A)) .
\enq

\mn
{\em Proof of Lemma~\ref{ML}.}
As already noted, Theorem~\ref{thmDT} gives Lemma~\ref{ML} when
$a_\f > M/(100D)$, so we assume this is not the case.

We assume (w.l.o.g.) that $x=n$ is a maximum degree
vertex of $\f$ and set
$\A  = \f\nex$ and
\[
\B
%= (\K_x\sm \f)^c
= \{V\sm T:T\in \K_x\sm \f\}
\]
%(the complements of sets on $x$ and not in $\f$, so
(so $|\A|=|\B|=a_\f$).

As above we take $\mm = n-1$.
The rest of our discussion takes place in the universe $V\sm \{x\} = [\mm]$.
%For the rest of the discussion we work with subsets of $[n]\sm \{x\} = [2k]$.
We use $\gG_l $ for $\C{[\mm]}{l}$---thus $\A\sub \gG_k$ (our earlier $V_{m,k}$)
%(used only with $l\in\{k,k+1\}$);
and $\B\sub \gG_{k+c}$---and set $\bar{\A}=\gG_k\sm \A$ and $\bar{\B}=\gG_{k+c}\sm \B$.

\medskip
For $\sss\sub \gG_k$ and $\T\sub \gG_{k+c}$, set
\[
\gL(\sss,\T) =|\{(A,B)\in \sss\times\T: A\sub B\}|.
\]
Notice that
\beq{gLe}
e(\f)    = \gL(\A,\bar{\B}) +e(\A) \geq \gL(\A,\bar{\B}).
\enq

We next observe that 
lower bounds on the $\gb$'s
imply lower bounds on the quantities $\gL(\A,\bar{\B})$:
\begin{prop}\label{betaProp}
For any $\f\in\cee$,
\beq{mainobs}
\gL(\A,\bar{\B}) \geq
\max\left\{(2\Cc{k}{c})^{-1}\gb_c(\A),\tfrac{1}{2ck}\Cc{k+c-2}{c-1}\gb_1(\A)
\right\}
 - \Cc{k+c-1}{c-1}|\A|/2.
\enq
\end{prop}

\nin
Of course in view of \eqref{gLe} this gives the same lower bound on $e(\f)$.

\mn
{\em Proof.}
The combination of
\[
\gL(\A,\B) +\gL(\A,\bar{\B}) = \gL(\A,\gG_{k+c}) =
\Cc{k+c-1}{c}|\A|
\]
and
\[
\gL(\A,\B) +\gL(\bar{\A},\B) =\gL(\gG_k,\B) =
\Cc{k+c}{c}|\B| =\Cc{k+c}{c}|\A|
\]
gives
\beq{mo1}
\gL(\bar{\A},\B) =\gL(\A,\bar{\B})+\Cc{k+c-1}{c-1}|\A|.
\enq

For the second bound in \eqref{mainobs}
we work in the (``Johnson") graph $J_1(m,k)$.
Write $\Phi$ for the number of triples
$(A,B,A')\in \A\times \gG_{k+c}\times \bar{\A}$
with $A\sim_1 A'$ and $A\cup A'\sub B$.
Since each relevant pair $(A,A')$ admits exactly $\C{k+c-2}{c-1}$ choices of $B$,
we have
\beq{mo2}
\Phi = \gb_1(\A) \Cc{k+c-2}{c-1}.
\enq
On the other hand, for each of the above triples, either
$(A,B)$ is one of the pairs counted by $\gL(\A,\bar{\B})$
or
$(A',B)$ is one of the pairs counted by $\gL(\bar{\A},\B)$
(and not both).
In the first case the number of choices of $A'$ is at most the number
of
%($J_1(m,k)$-)
neighbors of $A$ contained in $B$,
namely $ck$, and similarly in the second case.  This with
\eqref{mo1} gives
\[
\Phi \leq (\gL(A,\bar{\B})+\gL(\bar{\A},\B))ck
 = ck (2\gL(A,\bar{\B}) + \Cc{k+c-1}{c-1}|\A|),
\]
and then combining with
\eqref{mo2} yields the stated bound.

The argument for the first bound 
is similar and we just indicate the changes.
We work in $J_c(m,k)$ and consider triples as above but with $A\sim_c A'$
(so $B=A\cup A'$).  The number of triples, which is now just
$\gb_c(\A)$, is bounded above by
\[
(\gL(A,\bar{\B})+\gL(\bar{\A},\B))\Cc{k}{c}
 = \Cc{k}{c}(2\gL(A,\bar{\B}) + \Cc{k+c-1}{c-1}|\A|)
\]
($\C{k}{c}$ being the number of neighbors---now in $J_c(m,k)$---of $A$ contained in $B$
when $A\in \gG_k$, $B\in \gG_{k+c}$ and $A\sub B$), and the desired
bound follows.
\qqed

Finally, combining \eqref{mainobs} with \eqref{gLe} and our earlier bounds on the $\gb$'s
(see \eqref{A-M}-\eqref{betalb}) yields (with $\gc$ as in \eqref{betalb})
%(after some minor calculations)
\beq{mainobs2}
e(\f) \geq
\tfrac{|\A|}{2}\Cc{k+c-2}{c-1}
\max
\{\left(1-\tfrac{m}{k}\mu_k(\A)\right),
\tfrac{\gc m}{ck}\log\tfrac{1}{\mu_k(\A)}- \tfrac{k+c-1}{k}
\}.
\enq
(Replacing $\gb_c(\A)$ in \eqref{mainobs} by the lower bound provided by
\eqref{A-M} and \eqref{lambda} gives
\begin{eqnarray*}
e(\f) &\geq&
(2\Cc{k}{c})^{-1}
\tfrac{m}{k}\Cc{k}{c} \Cc{m-k-1}{c-1}|\A|(1-\mu_k(\A))
 - \Cc{k+c-1}{c-1}\tfrac{|\A|}{2}  \\
&=& \tfrac{|\A|}{2}\left[
\tfrac{m}{k}\Cc{m-k-1}{c-1}(1-\mu_k(\A))
 - \Cc{k+c-1}{c-1} \right]\\
&=& \tfrac{|\A|}{2} \Cc{k+c-2}{c-1}\left[
1-\tfrac{m}{k}\mu_k(\A) \right],
\end{eqnarray*}
and replacing $\gb_1$ by the lower bound from \eqref{betalb} yields
\[ %\begin{eqnarray*}
e(\f) 
\geq
\tfrac{1}{2ck}\Cc{k+c-2}{c-1}
\gc \mm |\A|
\log\tfrac{1}{\mu_k(\A)}    %\log ( 1/\mu_k(\A)) 
 - \Cc{k+c-1}{c-1}\tfrac{|\A|}{2},
\]
which is easily seen to be equal to the second bound in \eqref{mainobs2}.)

\medskip
It only remains to observe that this does what we want, namely
that (for suitable $\vt$) the expression in \eqref{mainobs2} is at least as large
as the bound in \eqref{e(f)}, or, equivalently, that
the max in \eqref{mainobs2} is at least

\beq{target}
\tfrac{2(k+c-1)}{ck}\vt \log (N/a_\f) < (4/c)\vt \log(1/\mu_k(\A));
\enq
we assert that this is true provided $\vt < \gc/5$.
%with the relevant term from the max depending on whether 

If $\log(1/\mu_k(\A)) \leq c/\gc$, then the r.h.s. of \eqref{target} is less than
the first term in the max
(which is essentially 1 since $\mu_k(\A) = |\A|/N< |\A|/M $, which we are assuming
is less than $ .01D^{-1}$; see following Theorem~\ref{thmDT}).

If, on the other hand, $\log(1/\mu_k(\A)) > c/\gc$, then
the second term in the max is at least
\[
\tfrac{\gc}{c}\tfrac{2k+c-1}{k}\log\tfrac{1}{\mu_k(\A)}- 
\tfrac{\gc}{c}\tfrac{k+c-1}{k} \log\tfrac{1}{\mu_k(\A)} =
\tfrac{\gc}{c}\log\tfrac{1}{\mu_k(\A)},
%> (4/c)\vt \log(1/\mu_k(\A)).
\]
which is again greater than the r.h.s. of \eqref{target}.\qqqed

\medskip
For $n=2k+1$ we could avoid the machinery 
used above for intermediate values of $|\A|$
(namely \eqref{A-M}, \eqref{lambda} and the first bound in \eqref{mainobs})
by choosing $\gz$ in Theorem~\ref{thmDT} to handle 
$\log(1/\mu_k(\A)) \leq c/\gc$ %(equivalently $|\A|\geq e^{-c/\gc}N$)
(and adjusting $\vt$ accordingly).

%(One final remark:  our original, somewhat harder proof for $n=2k+1$---found before 
%we were aware of \cite{Das}---depended crucially on the main result of \cite{ODW}.)

\bn
Department of Mathematics\\
Rutgers University\\
Piscataway NJ 08854\\
prd41@math.rutgers.edu\\
jkahn@math.rutgers.edu

\end{document}